\documentclass[pdflatex,sn-mathphys-num]{sn-jnl}
\usepackage{graphicx}%
\usepackage{multirow}%
\usepackage{amsmath,amssymb,amsfonts}%
\usepackage{amsthm}%
\usepackage{mathrsfs}%
\usepackage[title]{appendix}%
\usepackage{xcolor}%
\usepackage{textcomp}%
\usepackage{manyfoot}%
\usepackage{booktabs}%
\usepackage{algorithm}%
\usepackage{algorithmicx}%
\usepackage{algpseudocode}%
\usepackage{listings}%

\usepackage{graphicx} % Required for inserting images
\usepackage[utf8]{inputenc}
\usepackage{geometry,amsmath,amsthm,amssymb}
\usepackage{amsthm,pdfsync,enumerate,subfigure}
\usepackage[framemethod=tikz]{mdframed}
\usepackage{tikz-cd}
\hypersetup{allcolors=blue}
\numberwithin{equation}{section}
\newcommand{\I}{\mathrm{i}}

\newcommand{\lb}{\left(}

\newcommand{\rb}{\right)}
\newcommand{\PD}{\partial}

\newcommand{\Beq}{\begin{equation}}
\newcommand{\Eeq}{\end{equation}}
\newcommand{\beq}{\begin{equation*}}
\newcommand{\eeq}{\end{equation*}}
\newcommand{\bal}{\begin{align}}
\newcommand{\eal}{\end{align}}

\newcommand{\g}{\gamma}

\usepackage{mathtools}

{\definecolor{pink}{rgb}{1.0,.33.,64}

\newcommand{\B}{\beta}
\newcommand{\bp}{\begin{prob}}
	\newcommand{\ep}{\end{prob}}
\newcommand{\bpr}{\begin{proof}}
	\newcommand{\epr}{\end{proof}}
\renewcommand{\o}{\omega}

\newcommand{\st}{\,:\,}

\newcommand{\bel}[1]{\begin{equation}\label{#1}}
\newcommand{\ee}{\end{equation}}

\theoremstyle{definition}

\newcommand{\D}{\mathrm{d}}

\newcommand{\A}{\alpha}
\newcommand{\vp}{\varphi}

\usepackage{pdfsync,verbatim,epstopdf,enumerate}

\newcommand{\wt}{\widetilde}

\renewcommand{\o}{\omega}
\begin{document}

\title[Microlocal Analysis of 3D SAR Imaging]{Artefact Analysis of Multistatic Three-Dimensional SAR Imaging with Two Linear Trajectories} 


\author*[1]{\fnm{Daniel} \sur{Andre}}\email{d.andre@cranfield.ac.uk}

\author[2]{\fnm{Venkateswaran P.\ } \sur{Krishnan}}\email{vkrishnan@tifrbng.res.in}

\author[3]{\fnm{Clifford} \sur{Nolan}}\email{clifford.nolan@ul.ie}
\equalcont{All authors contributed equally to this work}
\affil*[1]{\orgdiv{Centre for Electronic Warfare, Information and Cyber}, \orgname{Cranfield University}, \orgaddress{\street{Defence Academy of the United Kingdom}, \city{Shrivenham}, \postcode{SN6 8LA}, \country{UK}}}

\affil[2]{\orgdiv{TIFR Centre for Applicable Mathematics}, \orgaddress{\street{PO Box 6503 GKVK PO}, \city{Bangalore}, \postcode{560065}, \state{Karnataka}, \country{India}}}

\affil[3]{\orgdiv{Department of Mathematics \& Statistics}, \orgname{University of Limerick}, \orgaddress{\street{Plassey Park Road}, \city{Castletroy}, \postcode{V94 T9PX}, \state{Limerick}, \country{Ireland}}}

\abstract
{In recent years, radar technology has seen much improvement, making multistatic Synthetic Aperture Radar (SAR) sensing a realistic possibility, for example in satellite constellations or unmanned aircraft systems. With such systems, there then comes the requirement to investigate useful multistatic SAR geometries. We provide a microlocal analysis of a multistatic data acquisition geometry for three-dimensional SAR imaging, where the transmitter and receiver move independently along two straight lines covering a region of interest (ROI). Using microlocal techniques, we show how the data acquisition geometry influences whether artefacts are likely to be present in the resultant image and how they can be avoided. Our main contribution is to provide a time-gating condition on the data which ensures that artefacts are not present in the ROI. As an independent verification, we provide several numerical simulations which follow a time-independent formulation.

% We provide a microlocal analysis of a multistatic data acquisition geometry in a three-dimensional SAR imaging set-up where the transmitter and receiver move independently along straight lines covering a region of interest (ROI) in three-dimensional space. Using microlocal techniques, we show how the data acquisition geometry influences whether artefacts are likely to be present in the image and how they can be avoided. Our main contribution is time-gating of the data in such a way that ensures that artefacts are not present in the ROI. As an independent verification, we provide several numerical simulations, following a time-independent formulation, to validate our results.  

%In recent years, a lot of work has focused on understanding how to efficiently model radar data and correspondingly, how to extract a reliable image of a region of interest (ROI) from this data.  

%Such artefacts would suggest that a fictitious object is present in the ROI when in fact that is not the case. 
}

%This situation often comes about due to symmetries within the data that specifically relate to the data acquistion geometry.  }   
\maketitle
\section{Introduction}
In recent years, much research has been conducted on bistatic radar imaging \cite{FundRadIm,YYC,ambartsoumian.felea.ea.18, felea.krishnan.ea.16, felea.gaburro.greenleaf.nolan.22}. The increased interest is partly due to improvements in technology which have made the proposition of multistatic radar sensing systems realisable. Microlocal analysis \cite{grigis.sjostrand, hormander.03, duistermaat.11, Treves} has proven to be a powerful tool for understanding how backprojected images of scenes are formed and how image artefacts can appear. Understanding how and where these artefacts may form is vital so that one may have the opportunity to carefully design how the Radar data is collected (the acquisition geometry) in order that such artefacts may be avoided in certain parts of the scene; the so-called regions of interest (ROI).

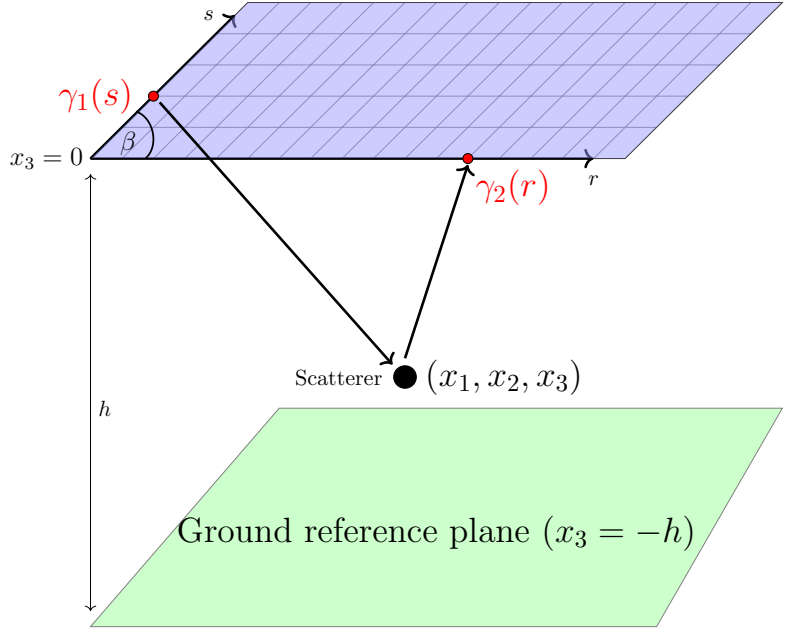
\begin{figure}[htbp]\label{Acquisition_Geometry_Figure}
  \centering
\resizebox{0.8\linewidth}{!}{  
\begin{tikzpicture}
  \draw[very thick,->] (3,10) -- (11,10) node[below, inner ysep=7pt] {$r$};
  \draw[very thick,->] (3,10) -- +(45:3.25cm) node[left, inner xsep=9pt]{$s$};
  \draw[shift={(3,10)},xslant=tan(45), step=.5, opacity=0.5] grid ++ (8.5,2.5);
  \draw[fill=blue!40!white, opacity=0.5] (3,10) -- (11.5,10) -- (14,12.5) -- (5.5,12.5) -- cycle;
  \draw[fill=green!40!white, opacity=0.5] (3,2.5) -- (12,2.5) -- (14,6) -- (6,6) -- cycle;
  \draw[<->] (3,9.75) -- (3,6) node[right]{$h$} -- (3,2.75);
  \draw[fill=red] (4,11) circle (0.08) node[left, inner xsep=9pt, font=\Large, color=red] {$\gamma_1(s)$};
  \draw[fill=red] (9,10) circle (0.08) node[below right, inner ysep=5pt, font=\Large, color=red] {$\gamma_2(r)$};
  \draw[fill=black] (8,6.5) circle (0.18) node[right, inner xsep=9pt, font=\Large, color=black] {$(x_1,x_2,x_3)$} node[inner xsep=10pt, left]{Scatterer};
  \draw[very thick, ->] (4.1,10.9) -- (7.8,6.7);
  \draw[very thick, ->] (8.0,6.8) -- (9,9.9);
  \draw (8.5,4) node[font=\Large]{\text{Ground reference plane ($x_3=-h$)}};
  \draw (2.3,10) node[font=\large]{\text{$x_3=0$}};
  \draw (3.6,10.25) node[font=\large]{\text{$\beta$}};
%  \draw[thick, black] (3.55,9.95) arc (23:41:1.35);
\node (1) at (3.75,9.9) {};
\node (2) at (3.6,10.8) {};
\draw  [thick, bend right=60] (1) to (2);
\end{tikzpicture}
}
\caption{Illustration of the proposed acquisition geometry, showing generic positions of the transmitter $\gamma_1(s)$, receiver $\gamma_2(r)$ and a scatterer at $x$ with coordinates $(x_1,x_2,x_3)$.}
\end{figure}

We consider multistatic (multiple bistatic) synthetic aperture radar (SAR) imaging that employs radio wave sources and receivers distributed along a pair of horizontal straight lines ($\gamma_1, \gamma_2$) which are at an angle $\beta$ to each other and these lines are parametrised by $s$ and $r$ respectively.  Thus, for $N$ transmitter positions and $M$ receiver positions, we have $MN$ bistatic radar pulse datasets. The SAR geometry is illustrated in Figure \ref{Acquisition_Geometry_Figure}. Such a multistatic SAR geometry would provide a virtual two-dimensional SAR aperture, so we expect, in principle, to be able to  reconstruct three-dimensional objects, as we have three degrees of freedom ($s,r$ - the slow time variables and  $t$ - the fast time variable). 

The ability to scan a two-dimensional virtual aperture, whilst constraining the transmitter and receiver (or transceiver) positions to lie only on two single lines could prove beneficial in practical scenarios where the region of space in between the lines may be inaccessible to sensors. It may also be the case that accurately arranging or scanning a transmitter or receiver purely along two straight lines may prove more practical to set-up when compared to accurately scanning a two-dimensional range of transmitter or receiver positions. 

It is noted that the work here, although directed to sensing with electromagnetic waves, is also relevant to synthetic aperture sensing with other wave modalities including sound waves, which has applications in underwater imaging, medical ultrasound, and seismic waves to image underground structures, used for geological exploration, for example.

There are at least three different ways that the multistatic collection we propose could be obtained in practice. To help in the understanding of the geometry, these are described in the following three cases.

Case A: One may have a set of $N$ fixed transmitters distributed over one line, and a set of $M$ fixed receivers distributed over the other. The transmitters would then transmit one pulse each in sequence, with the receivers providing $M$ recordings for each pulse;

Case B: As for Case A, but with the $N$ multiple transmitters may be substituted by a single transmitter that flies out along its line, transmitting $N$ pulses sequentially along its aperture, whilst the $M$ receivers record $M$ signals for each of the transmissions;

Case C: As for Case B, however the multiple receivers are substituted by a single receiver that flies along its aperture line multiple times. Here, one would arrange the collection so that for each of the $N$ transmitter positions, the transmitter generates a signal $M$ times, whilst the receiver does a full aperture sweep along its line to collecting the $M$ scattered transmitter pulses. For this case, only a single transmitter and a single receiver are required. 

%\tpk{[ Note: Are we aware of similar approaches being described in the literature? In SONAR and GPR perhaps? What can we say about the novelty here? Mention / reference "walk-away" paper? Cliff has added a paragraph about the Walkaway paper directly below this comment.]}

One of the acquisition geometries considered here was also considered in \cite{felea.gaburro.greenleaf.nolan.22}, where the potential for extra artefacts to occur in the image was identified as a possibility.  These artefacts may arise due to the breakdown of the local condition in the so-called B\"olker condition. More specifically, the B\"olker  condition requires a certain projection map from the wavefront relation of the scattering operator into the cotangent space of the data manifold to be an embedding; please see the discussion at the end of Section \ref{S2} for the terminology used here. The latter paper showed that the projection may fail to be an immersion (the local condition).  This paper revisits the B\"olker condition and shows that by designing the acquisition geometry appropriately (the extent and location of the source and receiver line segments), then it is possible to ensure that the B\"olker condition is satisfied so that we can then be confident regarding the absence of the artefacts mentioned in \cite{felea.gaburro.greenleaf.nolan.22}. In particular, in this paper, we provide a detailed calculation that shows how one may guarantee that the projection map is injective, putting wavefront set elements of the visible scatterers in a one-to-one correspondence with wavefront set elements of the data. We also consider in this paper the situation where the acquistion geometry has the sources and receivers arranged in two non-orthogonal straight line segments, thereby extending previous results, while at the same time showing that the reconstructions are stable with respect to the acquistion geometry.

We wish to analyse the image that could be obtained of scatterers, which in practice could include aircraft, buildings, ground terrain, etc., in the ROI (scene).  This image is obtained by backprojecting  the scattered wave field recorded at the receivers, which we call the data henceforth.  We will presently explain what is meant by the process of backprojection but for now, suffice it to say that it is a concept that is omnipresent in almost all imaging methods.  For example, in seismic imaging, we see the concept of ``migrating" the seismograms to form an image; see \cite{symes.09} for example. In the radar community, the backprojection image formation is utilized, and has been implemented on the simulated results described in Section \ref{S5} \cite{JE}.

We are particularly interested in investigating whether the proposed SAR geometry would give rise to artefacts (false scatterers) in the image.  When such artefacts occur, their contribution to the data is indistinguishable from that of a scatterer which is truly present in the scene, and hence it is imperative to determine their presence.

Our approach in this investigation is both theoretical through the use of microlocal analysis and also numerical to validate the theoretical predictions.  Using microlocal analysis, we study the propagation of the non-smooth (singular) part of the radio waves (which concentrates at the wavefronts) and how the corresponding singular component of the scattered waves can be used to infer the non-smooth component of the wave velocity field.  The non-smooth component of the velocity field is a function of the refractive indices of the materials (such as air, building materials, terrain composition, etc.) encountered by the propagating radio waves.  So ultimately, the image obtained is a representation of the estimated material composition of the scene.

Let us assume that before and after scattering, the radio waves are travelling in air with a constant speed of propagation $c_0$.  When the waves enter other materials where the speed of propagation is $c$, then we are interested in estimating the non-smooth part of the deviation $c_0-c$ from the background speed, which then provides us with an image of the scene. It turns out to be more mathematically convenient to instead work with the difference $c_0^{-2}-c^{-2}$, which we call the reflectivity function, $V$,  and the goal of the paper is to estimate the singular component of $V$ from the scattered waves.

We will show how the scattered wavefield can be represented by a so-called Forward Operator, $F$, so that $F:V\to d$  provides us with a model of the data $d$.   We will show that, under standard assumptions, the operator $F$ is a Fourier Integral Operator (FIO). In essence, the image of the scene is normally obtained by applying a weighted adjoint operator $F^*$ of $F$ to the data $d$.  The application of $F^*$ is essentially known as backprojection.  A rough interpretation of backprojection is a weighted summation of the various parts of the data that contribute to a particular point in the scene, where we wish to construct $V$. This is the same principle that many imaging methods are based upon, and the principle behind Global Navigation Satellite System positioning, commonly referred to as GPS.

From previous experience (see \cite{felea.krishnan.ea.16, gaburro.nolan.ea.07}, \cite{ambartsoumian.felea.ea.18}, for example) we know the operator $F$ is heavily influenced by the acquisition geometry, i.e. the location of the sources, receivers and recording interval used.  The acquisition geometry we use in this study employs sources along a line segment and receivers are along another line segment at the same altitude as the sources, with the two line segments intersecting at an angle $\beta$ to each other.  The reason for considering this set-up is that, intuitively, if the sources and receivers are collinear ($\B=0$ or $\B=\pi$ in our set up), a single scattering location would give rise to a whole ring of artefacts by symmetry considerations. Although this may not be immediately obvious, it will become clear after the analysis presented, with the reason being due to symmetries in the acquisition geometry.  By inclining the line segments to each other, we break this cylindrical symmetry and thereby avoid such circular artefacts.  It will turn out that it is possible to do so with a carefully arranged acquisition geometry, combined with certain filters applied to the data $d$ as a pre-processing step. However, constraining the sources and receivers to a plane will still give rise to a artefacts which are mirror images of true scatters, reflected across the plane generated by the lines of sources and receivers.

Because the image we obtain is essentially $F^*d=F^*FV$, the normal operator $F^*F$ becomes important to study. The sources and receivers lie at a fixed altitude and provided that we know that the scatterers recorded in the data are either below or above this height then our results may be summarized as follows. For a given aspect ratio between the line of sources and receivers, we identify a temporal window where no artefacts will be introduced into the scene.  Conversely, for a given region of interest, we may choose an appropriate aspect ratio and temporal window in order to avoid such artefacts.  If desirable, one may design a filter to pre-process the data so that it is unnecessary to place any further restrictions on the aspect ratio or temporal window wherein the data is recorded. 

In Section \ref{S2} the forward model is described. Section \ref{S3} provides the analysis of perpendicular transmitter and receiver axes, whereas Section \ref{S4} covers the non-perpendicular case generalization. Supporting numerical simulation results for four multistatic geometry scenarios are provided in Section \ref{S5}. The simulations follow a time-independent formulation different to that of the previous sections, but closely aligned with common radar laboratory measurement setups \cite{WAF, FundRadIm}. This also serves as independent verification of the methods presented in Sections \ref{S3} and \ref{S4}. Conclusions are provided in Section \ref{S6}.

\section{The Forward Model}\label{S2}
In this section, we recall a standard mathematical model, which uses geometrical optics to describe (i) the emission of radio waves from the transmitters, (ii) an approximation for the wave field that has scattered from inhomogeneities in the environment and (iii) the measurements of the resulting scattered waves at receivers.  This whole process is captured by the so-called scattering operator, described below.

Let us initially suppose that the transmitter and receiver move along a pair of orthogonal lines at the same altitude.  We parametrise these locations, using functions $\gamma_1:(s_{\mathrm{min}},s_{\mathrm{max}})\to \mathbb R^3$ and $\gamma_2:(r_{\mathrm{min}},r_{\mathrm{max}})\to \mathbb R^3$ with
\begin{eqnarray}
    \gamma_1(s) = (0,2s,0)\\ \gamma_2(r) = (2r,0,0).
\end{eqnarray}
We have set the altitude to be zero for mathematical convenience, meaning that the ground is located at $z=-h$, where $h>0$ represents the height at which the transmitters and receivers are located.

As each component of the electric field emitted by the transmitters satisfies a scalar wave equation, we will model the radio waves using a scalar wave equation. Such models are common in the literature; see for example \cite{ambartsoumian.felea.ea.18}, which we adapted here to obtain following approximate formula for $u(s,r,t)$, for waves emitted from the source location $\gamma_1(s)$, scattered and then measured at receiver location $\gamma_2(r)$ at time $t$: 
\begin{equation}\label{fm}
    u(s,r,t) = \int e^{i\Phi(s,r,t, x,\omega)}a(s,r,t,x,\omega)V(x)\ d\omega dx
\end{equation}
where the amplitude $a(s,r,t,x,\omega)$ is a symbol that captures the geometrical spreading of the energy along the incident and scattered wavefronts and the phase function 
\begin{equation} \label{Born}
    \Phi(s,r,t,x,\omega) = \omega(t-|\gamma_1(s)-x|/c_0-|x-\gamma_2(r)|/c_0)
\end{equation}
comes from the geometrical optics model for high frequency ($\omega \gg 1$) wave propagation. This amplitude also smoothly cuts off the data in order to avoid artefacts in the backprojected image later. The function (or distribution) 
\begin{equation}
  V(x)=\frac{1}{c_0^2}-\frac{1}{c^2(x)}  
\end{equation}
is associated to the deviation of the wave speed $c$ from the background speed $c_0$ in air.  Such deviations are associated to changes in material properties, where one region of space is occupied by a (possibly inhomogeneous) material, other than air. It turns out that instead of tracking the difference $c-c_0$, it is more convenient to track the difference of their squared reciprocals.  We observe that our sources are idealized as point sources (see \cite{ambartsoumian.felea.ea.18}) in space and time. The approximation \eqref{Born} is based on the Born approximation for single-scattered waves. We also note that \eqref{fm} is an {\em oscillatory integral} and as such may represent both a function as well as a  distribution.  

We consider measurements of $u(s,r,t)$, taken over the acquisition manifold 
\begin{equation}
    \mathcal M = \left\{\ (s,r,t)\ |\ s\in (s_{\rm min},s_{\rm max}), \ r\in (r_{\rm min},r_{\rm max}), \ t\in (0,T)\ \right\},
\end{equation}
where $T>0$ is the duration of the measurements. Since the amplitude $a$ and phase $\Phi$ in (\ref{fm}) satisfy the standard conditions \cite{duistermaat.11} for a symbol and non-degenerate phase function, we have that the following forward model operator 
\begin{equation}
    F:V\mapsto u|_{(s,r,t)\in \mathcal M}
\end{equation}
is a Fourier Integral Operator (FIO) \cite{duistermaat.11}.  We will exploit the well-known properties of FIOs, to analyse the image obtained by backprojection methods later in this paper. The canonical relation of $F$ is given by
 \[
C=\Bigg{\{}(s,t,\PD_{s}\Phi,\PD_{t}\Phi; x,-\PD_{x}\Phi)\st
\PD_{\o} \Phi=0\Bigg{\}}\] where $\PD_{s}$ is the differential in the $s$
variables, etc. \cite{Treves}.  The formal $L^2$
adjoint operator, $F^*$, is also an FIO associated with the canonical relation $C^t$ (which is $C$
with $T^*X$ and $T^*Y$ coordinates flipped, if $C=\{x, \xi; y , \eta\}$ then $C^t=\{ y, \eta; x, \xi \}$).

As already mentioned, a standard reconstruction method in imaging is to study the normal operator $F^{*}F$ and study what this normal operator does to the singularities of $V$. In order to perform this analysis, we study the composition of the canonical relations, $C^t\circ C$. This is given as follows: 
\Beq\label{CR-composition}
C^{t}\circ C=\{(x,\xi, \wt{x},\wt{\xi}\,): y=\wt{y}, \eta=\wt{\eta}\},
\Eeq
with $C^{t}=\{(x,\xi; \wt{y},\wt{\eta})\}$ and $C=\{(y, \eta; \wt{x},\wt{\xi}\,)\}$. 
%\tblue{[Possible confusion here with $x_1$ being regarded as the first cooredinate of $x$, etc?]}
\begin{figure}[h]\label{Fig2}
\centering
\begin{tikzcd}
 & C\arrow[ld,"\pi_{L}"'] \arrow[rd,"\pi_R"]\\
 T^{*}Y& & T^{*}X
\end{tikzcd}
\caption{Canonical projections}
\end{figure}
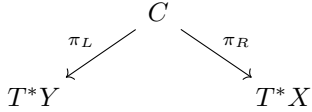
The composition of two FIOs is not always an FIO; unlike the situation for the case of pseudodifferential operators ($\Psi$DOs), which are special subclass of FIOs. In such instances, one is interested in analyzing the left and right projection maps $\pi_{L}$ and $\pi_{R}$, see Fig. \ref{Fig2}, to understand whether the composition of the two FIOs $F$ and $F^{*}$ introduces artefacts in the reconstructed images. One instance where the composition becomes an $\Psi$DO is the case when the left projection map $\pi_{L}$ is an embedding, that is, it is a proper injective immersion. This is the so-called B\"olker condition. In this special situation, we can be confident that no artefacts are introduced in the backprojected image.

\section{Analysis in the case of transmitter and receiver trajectories in perpendicular directions}\label{S3}

In this section, we will analyze the composition of the forward scattering operator $F$ with this formal $L^2$ adjoint $F^{*}$. Since we are interested in how these operators propagate singularities from the scene to the data (through the forward operator) and back to the scene (through the adjoint operator), we will analyze the composition at the level of composition of their corresponding canonical relations.

For simplicity, we let $h=0$ below. The case of a positive fixed $h>0$ can be easily handled by a translation.
Recall that the transmitter locations are denoted by $(0,2s,0)$ and the receiver locations by $(2r,0,0)$ for $s,r\geq 0$. Based on the discussion in the previous paragraph, we consider the following three equations based on the composition of the canonical relations $C^t$ corresponding to $F^{*}$ and $C$ corresponding to $F$; see \eqref{CR-composition}.
\begin{align}
   \label{Eq1} & \sqrt{(x_1-2r)^2+x_2^2 + x_3^2} + \sqrt{x_1^2+(x_2-2s)^2+x_3^2}\\
   \notag &\hspace{0.5in}=\sqrt{(y_1-2r)^2+y_2^2 + y_3^2} + \sqrt{y_1^2+(y_2-2s)^2+y_3^2},
   \end{align}
   \begin{align}
    \label{Eq2}&\frac{x_1-2r}{\sqrt{(x_1-2r)^2+x_2^2 + x_3^2}}=\frac{y_1-2r}{\sqrt{(y_1-2r)^2+y_2^2 + y_3^2}},\\
    \label{Eq3} & \frac{x_2-2s}{\sqrt{x_1^2+(x_2-2s)^2 + x_3^2}}=\frac{y_2-2s}{\sqrt{y_1^2+(y_2-2s)^2 + y_3^2}}.
\end{align}
Our goal is to show that the only intersection points of these three surfaces are the points 
\begin{align}\label{Eq4}
    ((x_1,x_2,x_3)=(y_1,y_2,y_3), \quad ((x_1,x_2,x_3)=(y_1,y_2,-y_3).
    \end{align}
    Note that by straightforward verification, \eqref{Eq4} satisfies the three equations \eqref{Eq1} - \eqref{Eq3}.

We will attempt to show this by rotating the coordinate system with $x_3$ fixed and using prolate spheroidal coordinates. %Note that the center of 
%We consider the  ellipsoid defined by 
%\begin{align}\label{Eq5}
%\sqrt{(x_1-2r)^2+x_2^2 + x_3^2} + \sqrt{x_1^2+(x_2-2s)^2+x_3^2}=t
%\end{align}
%is $C=(r,s,0)$. We will rotate coordinates with $x_3$ fixed, that is, in the $x_1-x_2$ plane anticlockwise through an angle $\A$.

%We first translate the ellipsoid so that the center of the ellipsoid is the origin and then we rotate the axes of the ellipsoid so that the rotated axes are aligned along the standard axes. 
Let us denote $\tan \A=\frac{s}{r}$, $D=\sqrt{r^2+s^2}$ 
and define the following change of coordinates: 
\Beq\label{Eq6}
\begin{aligned}
    &\wt{x}_1= \cos \A(x_1-r) -\sin \A (x_2-s),\\
  &\wt{x}_2=\sin \A (x_1-r)+\cos \A (x_2-s),\\
    &\wt{x}_3= x_3.
\end{aligned}
\Eeq
In the coordinate system $(\wt{x}_1,\wt{x}_2, \wt{x}_3)$, with $\wt{y}_1, \wt{y}_2, \wt{y}_3$ defined similarly, the ellipsoid defined by \eqref{Eq1} is given by 
\begin{align}
    \label{Eq9}&\sqrt{(\wt{x}_1 -D)^2+\wt{x}_2^2+ \wt{x}_3^2} + \sqrt{(\wt{x}_1+ D)^2+\wt{x}_2^2+ \wt{x}_3^2}\\
    \notag &\hspace{0.5in}=\sqrt{(\wt{y}_1 -D)^2+\wt{y}_2^2+ \wt{y}_3^2} + \sqrt{(\wt{y}_1+ D)^2+\wt{y}_2^2+ \wt{y}_3^2}.
\end{align}

Similarly \eqref{Eq2} transforms to
\begin{align}
    \label{Eq15}
    \frac{\cos \A \wt{x}_1 +\sin \A\wt{x}_2-r}{\sqrt{(\wt{x}_1 -D)^2+\wt{x}_2^2+ \wt{x}_3^2}}=\frac{\cos \A\wt{y}_1 +\sin \A\wt{y}_2-r}{\sqrt{(\wt{y}_1 -D)^2+\wt{y}_2^2+ \wt{y}_3^2}}.
\end{align}
%In \eqref{Eq15}, the expressions for $(\wt{y}_1, \wt{y}_2, \wt{y}_3)$ are defined analogous to the case of $(\wt{x}_1, \wt{x}_2, \wt{x}_3)$. Let us check \eqref{Eq15}. We substitute the expression for $(\wt{x}_1, \wt{x}_2, \wt{x}_3)$ from \eqref{Eq6}, \eqref{Eq7} and \eqref{Eq8}. We get for the left hand side of \eqref{Eq15} to be 
%\[
%\frac{x_2-2r}{\sqrt{(x_1-2r)^2+x_2^2+x_3^2}}.
%\]
%Analogous expression holds for $(\wt{y}_1,\wt{y}_2, \wt{y}_3)$. Analogous expression holds for $(\wt{y}_1,\wt{y}_2, \wt{y}_3)$.
and  \eqref{Eq3} transforms to 
\begin{align}
    \label{Eq16}\frac{-\sin \A \wt{x}_1 +\cos \A\wt{x}_2-s}{\sqrt{(\wt{x}_1+D)^2+\wt{x}_2^2+ \wt{x}_3^2}}=\frac{-\sin \A \wt{y}_1 +\cos \A\wt{y}_2-s}{\sqrt{(\wt{y}_1+D)^2+\wt{y}_2^2+ \wt{y}_3^2}}.
\end{align}
All these are standard verifications. 

Let us now define 
\Beq\label{Eq20}
\begin{aligned}
    &\wt{x}_1=D\cosh \rho \cos \phi\\
   &\wt{x}_2=D\sinh\rho \sin \phi \cos \theta\\
    &\wt{x}_3=D\sinh\rho \sin \phi \sin\theta,
\end{aligned}
\Eeq
with $\rho>0, \phi\in [0,\pi]$ and $\theta\in[0,2\pi]$. 

The points $(\wt{y}_1,\wt{y}_2,\wt{y}_3)$ will be denoted in the same way, but with a prime. Let us substitute this into \eqref{Eq9}, \eqref{Eq15} and \eqref{Eq16}. We then get the following. Note that in \eqref{Eq24} and \eqref{Eq25} below, we have used $s=D \sin \A$ and $r=D\cos \A$.
\begin{align}
  \label{Eq23} & \cosh \rho-\cos \phi + \cosh \rho+\cos \phi=\cosh \rho'-\cos \phi' + \cosh \rho'+\cos \phi',\\
   \label{Eq24}&\frac{\cos \A\lb \cosh \rho \cos \phi-1\rb+\sin \A \sinh\rho \sin \phi \cos \theta}{\cosh \rho -\cos \phi}\\&\notag \hspace{0.5in}=\frac{\cos \A\lb \cosh \rho' \cos \phi'-1\rb+\sin \A \sinh\rho' \sin \phi' \cos \theta'}{\cosh \rho' -\cos \phi'},\\
   \label{Eq25}&\frac{-\sin \A(\cosh\rho \cos \phi+1)+ \cos \A \sinh\rho \sin \phi \cos\theta}{\cosh \rho + \cos \phi}\\
   &\notag \hspace{0.5in}=\frac{-\sin \A(\cosh\rho' \cos \phi'+1)+ \cos \A \sinh\rho' \sin \phi' \cos\theta'}{\cosh \rho' + \cos \phi'}.
\end{align}
From \eqref{Eq23}, we have $\cosh\rho=\cosh\rho'$ and this implies $\rho=\rho'$. Next, we start simplifying \eqref{Eq24} and \eqref{Eq25}, taking into account the fact that $\rho=\rho'$.

From \eqref{Eq24}, we get, 
\begin{align}
    \notag &\cos \A \cosh^2\rho(\cos \phi-\cos\phi') + \sin \A \cosh\rho \sinh\rho(\sin \phi \cos\theta-\sin\phi' \cos \theta')\\
   \label{Eq26} &+\cos\A(\cos \phi'-\cos \phi) + \sin \A \sinh\rho(\sin\phi'\cos \theta' \cos \phi-\sin\phi\cos \phi' \cos \theta)=0.
\end{align}
This can be simplified to 
\begin{align}\label{Eq3.17}
    \notag &\sinh \rho \cos \A (\cos \phi -\cos \phi') + \sin \A \cosh\rho (\sin \phi \cos \theta - \sin \phi' \cos \theta')\\
    &+ \sin \A (\sin \phi' \cos\phi \cos \theta' - \sin \phi \cos \phi' \cos \theta)=0.
\end{align}
Similarly, from \eqref{Eq25}, we get, 
\begin{align}
    \notag &\sin \A\cosh^2 \rho(\cos\phi'-\cos\phi)+\cos \A \sinh\rho\cosh\rho(\sin \phi \cos \theta-\sin\phi'\cos \theta')\\
   \label{Eq27} &+\sin \A(\cos \phi-\cos \phi') + \cos \A \sinh\rho(\sin \phi \cos \theta \cos \phi'-\sin \phi'\cos \phi \cos \theta')=0.
\end{align}
Simplifying this further, we get, 
\[
\begin{aligned}
&-\sin \A (\cos \phi -\cos \phi')\sinh^2 \rho + \cos \A \sinh\rho\cosh\rho(\sin \phi \cos \theta-\sin\phi'\cos \theta')\\
&+\cos \A \sinh\rho(\sin \phi \cos \theta \cos \phi'-\sin \phi'\cos \phi \cos \theta')=0.
\end{aligned}
\]
Canceling out $\sinh \rho$, we then get, 
\Beq\label{Eq3.19}
\begin{aligned}
    & -\sin \A \sinh \rho (\cos \phi -\cos \phi') + \cos \A \cosh \rho(\sin \phi \cos \theta-\sin\phi'\cos \theta')\\
    &+\cos \A (\sin \phi \cos \theta \cos \phi'-\sin \phi'\cos \phi \cos \theta')=0.
\end{aligned}
\Eeq

For ease of notation, let us denote the following:
\Beq
\begin{aligned}
    & A=\cos \phi - \cos \phi',\\
    &B=\sin \phi \cos \theta - \sin \phi' \cos \theta',\\
    &C=\sin\phi'\cos \theta' \cos \phi-\sin\phi\cos \phi' \cos \theta
\end{aligned}
\Eeq
Then with this, \eqref{Eq3.17} and \eqref{Eq3.19} transform to 
\begin{align}
    \label{Eq3.21}& \cos \A \sinh \rho A + \sin \A \cosh \rho B +\sin \A C=0,\\
    \label{Eq3.22}& -\sin \A \sinh \rho A + \cos \A \cosh \rho B - \cos \A  C=0.
    \end{align}
  %  Cancelling out the common term of $\sinh \rho$ from the two equations, we get, 
  %  \begin{align}
  %     \label{Eq33A}& \cos \A \sinh \rho A + \sin \A \cosh \rho  B +\sin \A C=0,\\
  %  \label{Eq34A}& -\sin \A \sinh \rho A + \cos \A  \cosh \rho B - \cos \A  C=0. 
  %  \end{align}
    \begin{comment}
Let us multiply \eqref{Eq33A} by $\sin \A$ and \eqref{Eq34A} by $\cos \A$ and add. After simplifying, we get,
\begin{align}\label{Eq35}
    \cosh\rho B-\cos 2\A C =0.
\end{align}
%\begin{align}
 %   \cosh \rho \sinh \rho (\sin \phi \cos \theta -\sin \phi' \cos \theta')+ \sinh \rho(\sin \phi \cos \theta \cos \phi'-\sin \phi' \cos \phi \cos \theta')=0.
%\end{align}
%Grouping terms, we get, 
%\begin{align}\label{Eq29}
  %  \sin \phi \cos \theta(\cosh \rho + \cos \phi')-\sin \phi' \cos \theta'(\cosh \rho + \cos \phi)=0.
%\end{align}
%Let us rewrite this as follows: 
%\begin{align}\label{Eq30}
   % \sin \phi \cos \theta\cosh \rho-\sin \phi' \cos \theta'\cosh \rho=\sin \phi' \cos \theta' \cos \phi-\sin \phi \cos \theta \cos \phi'.
%\end{align}
Next, let us multiply \eqref{Eq33A} by $\cos \A$ and \eqref{Eq34A} by $\sin\A$ and subtract. We then get 
\begin{align}\label{Eq36}
    \sinh \rho A + \sin 2\A C =0.
\end{align}
\end{comment}
Next we write 
\begin{align}\label{Eq35}
C=\sin \phi'\cos \theta' A -B \cos \phi'.
\end{align}
Considering \eqref{Eq3.21}, \eqref{Eq3.22} and \eqref{Eq35}, we have a ``linear'' system of equations: 
\begin{align}
   \begin{pmatrix}
       \cos \A \sinh\rho & \sin \A \cosh \rho & \sin \A\\
       -\sin \A \sinh \rho & \cos \A \cosh \rho & -\cos \A\\
    \sin \phi' \cos \theta' & -\cos \phi' & -1\end{pmatrix}
       \begin{pmatrix}
       A\\    
       B\\
       C
   \end{pmatrix}=0.
\end{align}
The determinant $\Delta$ of the matrix on the left hand side is 
\begin{align}
    \notag \Delta &= \cos \A \sinh \rho(-\cos \A \cosh \rho -\cos \A \cos \phi')\\
    \notag &-\sin \A \cosh \rho(\sin \A \sinh \rho+\sin \phi' \cos \theta' \cos \A)\\
    &+\sin \A (\sin \A \cos \phi' \sinh \rho -\cos \A \sin \phi' \cos \theta' \cosh \rho).
\end{align}
Simplifying this, we get, (up to a negative sign)
\[
\Delta = \cosh \rho \sinh \rho +\sinh\rho \cos 2\A \cos \phi'+\cosh \rho \sin 2\A \sin \phi' \cos \theta'.
\]
Let us now derive conditions under which this determinant is non-vanishing. We have 
\begin{align*}
    \Delta &\geq \cosh \rho \sinh \rho -\sinh \rho -\cosh \rho\\
    &= \frac{1}{4}\lb e^{2\rho} - e^{-2\rho}\rb - \frac{1}{2}\lb e^{\rho}-e^{-\rho}\rb -\frac{1}{2}\lb e^{\rho} +e^{-\rho}\rb\\
    &=\frac{1}{4} \lb e^{2\rho}-e^{-2\rho}\rb -e^{\rho}\\
    &=\frac{e^{4\rho}-4e^{3\rho}-1}{4 e^{2\rho}}.
\end{align*}
We can choose $\rho\gg 1$ such that $e^{4\rho}-4e^{3\rho}-1>0$. For instance, a non-optimal bound can be obtained as follows: Since $\rho>0$, we have $1< e^{3\rho}$. Then 
\begin{align*}
    e^{4\rho} - 4e^{3\rho}-1\geq e^{4\rho}-5e^{3\rho}= e^{3\rho}\lb e^{\rho}-5\rb.
\end{align*}
For example, choosing  any $\rho> \ln 6$, we get that 
\[
e^{4\rho}-4e^{3\rho}-1>0.
\]
Based on this analysis, we conclude that for large enough $\rho$, the only additional singularities arise exactly as in \eqref{Eq4}. We also note that we are explicitly ensuring that we avoid the potential artefacts that were identified in the ``walkaway" acquisition geometry discussed in \cite[\S 6.2]{felea.gaburro.greenleaf.nolan.22}. This calculation also gives extra insight as to how the artefact identified in \cite[\S 6.2]{felea.gaburro.greenleaf.nolan.22} can arise. %\tpk{[should comment a little more explicitly on this.]} 

\section{Analysis in the case of transmitter and receiver trajectories in non-perpendicular directions}\label{S4}

In the interest of generality, we will go through an analysis corresponding to that conducted in \S \ref{S3} for the case when the transmitter and receiver trajectories are not necessarily perpendicular, but still following straight-line trajectories. While the calculations in the previous section will follow as a special case of what is done below, we prefer to keep these sections separate in the interest of readability.

As before, we are interested in the composition of the forward scattering operator with that of its adjoint at the level of composition of their corresponding canonical relations. 

%\subsection{Analysis of the forward operator}

The forward operator in the case of non-perpendicular trajectories for the transmitter and receiver is 
\Beq
\begin{aligned}
    F V(s,r,t)=\int e^{\I \o(t-\frac{1}{c_0}\lb \lvert \g_1(s)-x\rvert+\lvert x-\g_2(r)\rvert\rb} a(s,r,x,\o) V(x) \D x \D \o,
\end{aligned}
\Eeq
where $\g_1(s)=(2s\cot \B, 2s,0)$ and $\g_2(r)=(2r,0,0)$. We note that $\cot \B=0$, corresponding to when $\B=\pi/2$, covers the case considered in \S \ref{S3}.

The canonical relation of this operator is 
\Beq\label{Eqs4.2}
\begin{aligned}
    C=\lb (s, r, t, \PD_{s,r,t}\vp); (x,-\PD_{x}\vp): \PD_{\o}\vp=0\rb,
\end{aligned}
\Eeq
where 
\[
\vp = \vp(s,r,t,x,\o)= \o(t-\frac{1}{c_0}\lb |x-\g_1(s)| + |\g_2(r)-x|\rb).
\]
We note that $(s,r,x,\o)$ is a global parameterization for $C$.

The trajectories of the transmitter and receiver are $\g_{T}(s)=(2s\cot \B, 2s , 0)$ and $\g_R(s)=(2r,0,0)$. First, let us assume that $\cot \B>0$. Some changes are required in the calculations below for the case $\cot \B<0$.

With the trajectories as above, we consider the following three equations; see \eqref{CR-composition}:  
\begin{align}
   \label{Eq2.1} \notag &\sqrt{(x_1-2s\cot\B)^2+(x_2-2s)^2+x_3^2}+ \sqrt{(x_1-2r)^2+x_2^2+x_3^2}\\
    &\hspace{0.5in}=\sqrt{(y_1-2s\cot\B)^2+(y_2-2s)^2+y_3^2}+ \sqrt{(y_1-2r)^2+y_2^2+y_3^2},\\
   \label{Eq2.2} &\frac{\cot \B(x_1-2s \cot \B)+(x_2-2s)}{\sqrt{(x_1-2s\cot\B)^2+(x_2-2s)^2+x_3^2}}=\frac{\cot \B(y_1-2s \cot \B)+(y_2-2s)}{\sqrt{(y_1-2s\cot\B)^2+(y_2-2s)^2+y_3^2}},\\
    \label{Eq2.3}&\frac{x_1-2r}{\sqrt{(x_1-2r)^2+x_2^2+x_3^2}}=\frac{y_1-2r}{\sqrt{(y_1-2r)^2+y_2^2+y_3^2}}.
\end{align}
Let us consider a generic ellipsoid of the form
\begin{align}\label{Eq2.4}
\sqrt{(x_1-2s\cot\B)^2+(x_2-2s)^2+x_3^2}+ \sqrt{(x_1-2r)^2+x_2^2+x_3^2}=t.
\end{align}
The foci of this ellipsoid are $\g_T(s)=(2s\cot \B, 2s , 0)$ and $\g_R(s)=(2r,0,0)$. The centre of this ellipsoid is $(s\cos \B+r, s, 0)$. We apply a translation so that the ellipse is centred at the origin, which means the new foci are located at $(s\cot\beta-r,s,0)$ and $(r-s\cot\beta,-s,0)$. Next we rotate the ellipsoid by an angle $\A$ defined below, so that the ellipsoid has an axis aligned with the $x_1$-axis. We define
\begin{align}
\delta = \cot\beta, \ 
 \tan \A = \frac{s}{r-s\delta}.
\end{align}
The translation and rotation mentioned above is implemented via the following change of coordinates: 
\begin{align}
    &\wt{x}_1= \cos \A (x_1-r-s\cot \B)-\sin \A (x_2-s)\\
    &\wt{x}_2=\sin \A (x_1-r-s\cot \B)+\cos \A (x_2-s)\\
    &\wt{x}_3=x_3.
\end{align}
The ellipsoid defined by \eqref{Eq2.4} is transformed in the $\wt{x}_1, \wt{x}_2$ and $\wt{x}_3$ coordinates to
\begin{align}
    \sqrt{\lb \wt{x}_1 +D\rb ^2 +\wt{x}_2^2+\wt{x}_3^2} + \sqrt{\lb \wt{x}_1 -D\rb ^2 +\wt{x}_2^2+\wt{x}_3^2}=t
\end{align}
where 
\[
D=\sqrt{s^2+(\delta s-r)^2}.
\]
There follows by a straightforward calculation similar to that for the case of perpendicular trajectories in \S \ref{S3}. Hence \eqref{Eq2.1} transforms to 
\begin{equation}\label{Eq2.9}
\begin{aligned} 
&\sqrt{\lb \wt{x}_1 +D\rb ^2 +\wt{x}_2^2+\wt{x}_3^2} + \sqrt{\lb \wt{x}_1 -D\rb ^2 +\wt{x}_2^2+\wt{x}_3^2}\\
 &\hspace{0.5in}=\sqrt{\lb \wt{y}_1 +D\rb ^2 +\wt{y}_2^2+\wt{y}_3^2} + \sqrt{\lb \wt{y}_1 -D\rb ^2 +\wt{y}_2^2+\wt{y}_3^2}.
 \end{aligned}
\end{equation}

Next, \eqref{Eq2.2} transforms to 
\begin{align}\label{Eq2.10}
& \frac{\delta\Big{\{}\cos \A \wt{x}_1+\sin \A \wt{x}_2+D\cos\alpha\Big{\}}+\Big{\{}-\sin\alpha\wt{x}_1+\cos\alpha\wt{x}_2-D\sin\alpha\Big{\}}}{\sqrt{\lb \wt{x}_1+D\rb ^2 +\wt{x}_2^2+\wt{x}_3^2}} \notag \\
= & \frac{\delta\Big{\{}\cos \A \wt{y}_1+\sin \A \wt{y}_2+D\cos\alpha\Big{\}}+\Big{\{}-\sin\alpha\wt{y}_1+\cos\alpha\wt{y}_2-D\sin\alpha\Big{\}}}{\sqrt{\lb \wt{y}_1+D\rb ^2 +\wt{y}_2^2+\wt{y}_3^2}}.
\end{align}
Similarly \eqref{Eq2.3} transforms to 
\begin{equation}\label{Eq2.11}
\frac{ \cos \A \wt{x}_1+\sin \A \wt{x}_2-D\cos\alpha}{\sqrt{\lb \wt{x}_1-D\rb ^2 +\wt{x}_2^2+\wt{x}_3^2}}=\frac{ \cos \A \wt{y}_1+\sin \A \wt{y}_2-D\cos\alpha}{\sqrt{\lb \wt{y}_1-D\rb ^2 +\wt{y}_2^2+\wt{y}_3^2}}
\end{equation}
We consider now the equations \eqref{Eq2.9}, \eqref{Eq2.10}, \eqref{Eq2.11} and use prolate spheroidal coordinates to simplify them.

We let 
\begin{align}\label{Eq2.12}
\notag & \wt{x}_1=D \cosh \rho \cos \phi,\\
&\wt{x}_2=D \sinh \rho \sin \phi \cos \theta,\\
\notag & \wt{x}_3=D \sinh\rho \sin \phi \sin \theta.
\end{align}
We next substitute \eqref{Eq2.12} into \eqref{Eq2.9}, \eqref{Eq2.10} and \eqref{Eq2.11}. We start with \eqref{Eq2.9}. We get (denoting prime for the $\wt{y}$ variables) 
\begin{align}
    & D(\cosh\rho + \cos\phi) +  D(\cosh\rho-\cos\phi)=\\
    &=D(\cosh\rho' + \cos\phi') +  D(\cosh\rho'-\cos\phi'),
\end{align}
which gives 
\begin{align}
    \cosh \rho=\cosh \rho'
\end{align}
and this implies that $\rho =\rho'$. 

Next substituting \eqref{Eq2.12} into \eqref{Eq2.10}, we get, 
\Beq
\begin{aligned}
    \notag &\frac{\delta \Big{\{}D \cos \A +D \cos \A \cosh\rho \cos \phi+D\sin \A \sinh\rho \sin \phi \cos \theta\Big{\}}}{D(\cosh\rho + \cos\phi)}\\
    &\hspace{1in}+\frac{\Big{\{}-D \sin \A\cosh\rho \cos\phi +D\cos \A \sinh\rho \sin\phi \cos \theta-D \sin \A\Big{\}}}{D(\cosh\rho + \cos\phi)}
    \\&=\frac{\delta\Big{\{}D \cos \A +D \cos \A \cosh\rho \cos \phi'+D\sin \A \sinh\rho \sin \phi' \cos \theta'\Big{\}}}{D(\cosh\rho + \cos\phi')}\\
    &\hspace{1in}+\frac{\Big{\{}-D \sin \A\cosh\rho \cos\phi' +D\cos \A \sinh\rho \sin\phi' \cos \theta'-D \sin \A\Big{\}}}{D(\cosh\rho + \cos\phi')}.
\end{aligned}
\Eeq
Cross multiplying, we obtain, 
\Beq\label{Eq3.17A}
\begin{aligned}
     & D\delta \cosh\rho\cos \A  + D \delta \cosh^2 \rho \cos \A \cos \phi + D\delta \sin \A \sinh \rho \cosh \rho \sin \phi \cos \theta\\
     &-s \cosh \rho -D\sin \A \cosh^2\rho \cos \phi + D \cos \A \sinh \rho \cosh \rho \sin \phi \cos \theta\\
     &+ D\delta \cos \phi'\cos \A  + D \delta \cos \A  \cosh \rho \cos \phi \cos \phi' + D\delta \sin \A  \sinh \rho \sin \phi \cos \phi' \cos \theta \\
     &-s \cos \phi'-D \sin \A \cosh\rho \cos \phi \cos \phi' + D \cos \A \sinh\rho \sin \phi\cos \phi' \cos \theta\\
    = & D \delta \cosh\rho \cos \A  + D \delta\cosh^2 \rho \cos \A  \cos \phi' + D\delta \sin \A \sinh \rho \cosh \rho \sin \phi' \cos \theta'\\
     &-s \cosh \rho -D\sin \A \cosh^2\rho \cos \phi' + D \cos \A \sinh \rho \cosh \rho \sin \phi' \cos \theta'\\
     &+ D \delta \cos \A  \cos \phi + D \delta \cos \A  \cosh \rho \cos \phi' \cos \phi + D\delta \sin \A \sinh \rho \sin \phi'\cos \phi \cos \theta \\
     &-s \cos \phi-D \sin \A \cosh\rho \cos \phi \cos \phi' + D \cos \A \sinh\rho \sin \phi'\cos \phi \cos \theta' = 0.
\end{aligned}
\Eeq

%Cross multiplying, we obtain, 
%\begin{align}
 %   \notag & D \delta \cosh\rho \cos \A  + D \cosh^2 \rho \cos \A \cot \B \cos \phi + D\cot \B \sin \A \sinh \rho \cosh \rho \sin \phi \cos \theta\\
  %  \notag &-s \cosh \rho -D\sin \A \cosh^2\rho \cos \phi + D \cos \A \sinh \rho \cosh \rho \sin \phi \cos \theta\\
   % \notag &+ \cot \B \cos \phi'(r-s \cot \B) + D \cos \A \cot \B \cosh \rho \cos \phi \cos \phi' + D\sin \A \cot \B \sinh \rho \sin \phi \cos \phi' \cos \theta \\
   %\notag  &-s \cos \phi'-D \sin \A \cosh\rho \cos \phi \cos \phi' + D \cos \A \sinh\rho \sin \phi\cos \phi' \cos \theta\\
    %=\notag & \cosh\rho \cot \B(r-s \cot \B) + D \cosh^2 \rho \cos \A \cot \B \cos \phi' + D\cot \B \sin \A \sinh \rho \cosh \rho \sin \phi' \cos \theta'\\
    %\notag &-s \cosh \rho -D\sin \A \cosh^2\rho \cos \phi' + D \cos \A \sinh \rho \cosh \rho \sin \phi' \cos \theta'\\
    %\notag &+ \cot \B \cos \phi(r-s \cot \B) + D \cos \A \cot \B \cosh \rho \cos \phi' \cos \phi + D\sin \A \cot \B \sinh \rho \sin \phi'\cos \phi \cos \theta \\
   %\notag  &-s \cos \phi-D \sin \A \cosh\rho \cos \phi \cos \phi' + D \cos \A \sinh\rho \sin \phi'\cos \phi \cos \theta'.
%\end{align}
In order to simplify terms involving \eqref{Eq2.10} and \eqref{Eq2.11},  we let, as before, 
\Beq\label{Eq4.18}
\begin{aligned}
A = \cos \phi - \cos \phi'\\
B = \sin \phi \cos \theta - \sin \phi' \cos \theta'\\
C=\sin \phi' \cos \theta' \cos \phi - \sin \phi \cos \phi' \cos \theta
\end{aligned}
\Eeq
Also, as before,  $C$ can be written as 
\Beq\label{Eq4.19}
C=\sin \phi'\cos \theta' A -B \cos \phi'.
\Eeq

Simplifying \eqref{Eq3.17A} using \eqref{Eq4.18}, we get, 
%\begin{align}
%\notag & D \delta \cosh^2 \rho \cos \A  A + D \delta \sin \A \cosh \rho \sinh \rho B -D \sin \A \cosh^2 \rho A + D \cos \A \cosh \rho \sinh \rho B\\
%&-D\delta \cos \A A - D\delta \sin \A  \sinh \rho C + s A  - D \cos \A \sinh\rho C=0.
%\end{align}
%Using the fact that $r+s\cot \B= D \cos \A$ and $s=D \sin \A$, we get,
%\begin{align}
%\notag & D \cosh^2 \rho \cos \A \cot \B A + D \cot \B \sin \A \cosh \rho \sinh \rho B -D \sin \A \cosh^2 \rho A + D \cos \A \cosh \rho \sinh \rho B\\
%&-D \cot \B \cos \A A +2D \sin \A \cot^2 \B A - D\sin \A \cot \B \sinh \rho C + D \sin \A A  - D \cos \A \sinh\rho C=0.
%\end{align}
%Simplifying this, we have, 
%\begin{align}
%\notag &  \cosh^2 \rho \cos \A \cot \B A +  \cot \B \sin \A \cosh \rho \sinh \rho B - \sin \A \sinh^2 \rho A +  \cos \A \cosh \rho \sinh \rho B\\
%&-\cot \B \cos \A A +2 \sin \A \cot^2 \B A - \sin \A \cot \B \sinh \rho C   -  \cos \A \sinh\rho C=0.
%\end{align}
%Further simplifying, we get, 
\begin{align}\label{Eq2.23}
    A\sinh \rho (\delta \cos \A-\sin \A)+ B \cosh \rho (\cos \A + \delta \sin \A)-C(\cos \A + \delta \sin \A)=0.
\end{align}
%Note that when $\cot \B=0$, i.e, $\delta=0$,  this equation matches with \eqref{Eq34A}.

Next we substitute \eqref{Eq2.12} into \eqref{Eq2.11} to get 
\begin{equation}
\begin{aligned} 
&\frac{\cosh\rho\cos\alpha\cos\phi +\sinh\rho\sin\alpha\sin\phi\cos\theta-\cos\alpha}{(\cosh \rho -\cos \phi)}\\
&\hspace{0.5in}= \frac{\cosh\rho\cos\alpha\cos\phi' +\sinh\rho\sin\alpha\sin\phi'\cos\theta'-\cos\alpha}{(\cosh \rho -\cos \phi')}.
\end{aligned} 
\end{equation}
Cross multiplying, we obtain, 
\[
\begin{aligned}
    & \cosh^2\rho \cos\alpha\cos\phi +\cosh\rho\sinh\rho\sin\alpha\sin\phi\cos\theta - \sinh\rho\sin\alpha\cos\phi'\sin\phi\cos\theta 
    \\
    & -\cosh^2\rho\cos\alpha\cos\phi'  -\cosh\rho\sinh\rho\sin\alpha\sin\phi'\cos\theta' +\sinh\rho\sin\alpha\cos\phi\sin\phi'\cos\theta' \\
    & +\cos\phi'\cos\alpha-\cos\phi\cos\alpha = 0.
\end{aligned}
\]
Simplifying this, we obtain, 
\begin{align}
(\cosh^2 \rho-1) \cos \A A + \cosh \rho \sinh\rho \sin \A B +\sinh\rho \sin \A  C =0.
\end{align}
Therefore we have
\begin{align}\label{Eq2.28}
A\sinh \rho \cos \A + \cosh \rho  \sin \A B+  \sin \A C =0.
\end{align}
Note that this matches with \eqref{Eq3.21}.

We combine \eqref{Eq4.19}, \eqref{Eq2.23} and \eqref{Eq2.28} into a system: 
\Beq
\begin{aligned}
 &A\sinh \rho (\delta \cos \A-\sin \A)+ B \cosh \rho (\cos \A + \delta \sin \A)-C(\cos \A + \delta \sin \A)=0,\\
    & A\sinh \rho \cos \A + \cosh \rho  \sin \A B+  \sin \A C =0\\
    &A\sin\phi' \cos \theta' -B\cos \phi' -C=0.
\end{aligned}
\Eeq
Let us write this in matrix form: 
\Beq
\begin{pmatrix}
    \sinh\rho(\delta \cos \A-\sin \A) & \cosh \rho (\cos \A + \delta \sin \A) &-(\cos \A + \delta \sin \A)\\
    \sinh\rho \cos \A & \cosh \rho \sin \A &\sin \A \\
    \sin \phi' \cos \theta' &-\cos \phi' &-1
\end{pmatrix}\begin{pmatrix}
    A\\ B\\C
\end{pmatrix}=0.
\Eeq
We compute the determinant $\Delta$ of the matrix on the left hand side. % We get, expanding along the third row: 
%\begin{align*}
 %   \Delta&=-\lb \cosh \rho \sinh \rho\sin \A(\delta \cos \A -\sin \A)-\cosh\rho\sinh \rho \cos \A(\cos \A + \delta \sin \A)\rb \\
  %  &+\cos \phi' \lb \sinh \rho \sin \A(\delta \cos \A -\sin \A)+ \sinh \rho \cos \A (\cos \A +\delta \sin \A)\rb \\
  %  &+\sin \phi' \cos \theta'\lb \cosh \rho \sin \A(\cos \A + \delta \sin \A) + \cosh\rho \sin \A(\cos \A + \delta \sin \A)\rb.
%\end{align*}
We get 
\begin{align*}
    \Delta &= \cosh \rho \sinh \rho +\sinh \rho \cos \phi' \lb \delta \sin 2\A +\cos 2\A\rb \\
    &\hspace{0.5in}+ \cosh\rho \sin \phi' \cos \theta' \lb \sin 2\A + 2\delta \sin^2 \A\rb.
\end{align*}

Next, let us impose conditions on $\rho$ to make this determinant non-vanishing. %Using the lower bounds for the trigonometric functions appearing in $\Delta$ above, 
We have,
\Beq
\begin{aligned}
\Delta &\geq \cosh \rho \sinh \rho -(1+|\delta|)\sinh \rho  -(1+2|\delta|)\cosh\rho\\ 
&=\frac{e^{2\rho}-e^{-2\rho}}{4} -\frac{(1+|\delta|)}{2}(e^{\rho}-e^{-\rho})-\frac{1+2|\delta|}{2}(e^{\rho}+e^{-\rho})\\
&=\frac{e^{-2\rho}}{4}\lb e^{4\rho}-1 -2(1+|\delta|)(e^{3\rho}-e^{\rho})-2(1+2|\delta|)(e^{3\rho}+e^{\rho})\rb \\
&=\frac{e^{-2\rho}}{4}\lb e^{4\rho}-1-(2+2|\delta|+2+4|\delta|)e^{3\rho} +(2+2|\delta|-2-4|\delta|)e^{\rho}\rb\\
&=\frac{e^{-2\rho}}{4}\lb e^{4\rho}-1-(4+6|\delta|)e^{3\rho}-2|\delta| e^{\rho}\rb.
\end{aligned}
\Eeq
We choose $\rho\gg 1$, so that $1<e^{3\rho}$ and $e^{\rho}<e^{3\rho}$. With this, we get, 
\[
\Delta > \frac{e^{-2\rho}}{4}\lb e^{4\rho}-(5+8|\delta|)e^{3\rho}\rb\\
=\frac{e^{\rho}}{4}\lb e^{\rho}-(5+8|\delta|)\rb.
\]
Given a fixed $0\leq |\delta|<\infty$, we can choose $\rho$ large enough so that $\Delta>0$. This would mean $A=B=C=0$. Beginning with $A=0$, we get $\cos \phi =\cos \phi'$. Since $0\leq \phi,\phi'\leq \pi$, we have that $\phi'=\phi$. Hence $\sin \phi = \sin \phi'$, which then gives that $\cos \theta =\cos \theta'$. We then have $\theta=\theta'$ or $\theta'= 2\pi -\theta$.

As with the analysis of the perpendicular trajectories case in \S \ref{S3}, we also conclude for the non-perpendicular case, that for large enough $\rho$, the only additional singularities arise exactly as in \eqref{Eq4}. We additionally remark that the calculations in this section generalize the analysis made for the walkway geometry in \cite[\S 6.2]{felea.gaburro.greenleaf.nolan.22} which only considered the perpendicular trajectories case.

\section{Numerical Simulations}\label{S5}

To demonstrate the proposed multistatic SAR aperture approach, a free-space simulation was conducted for a specific scenario, involving isotropic point scatterers, resulting in three-dimensional SAR imagery. The approach undertaken comprised bistatic and three-dimensional generalization of the time-independent approach described by Gorham et. al. \cite{GM}. The simulation approach aligns closely with the radar measurement approach utilized in the Ground-Based SAR laboratory based at Cranfield University \cite{WAF}\cite{FundRadIm}, which utilizes a stepped-frequency chirp waveform created by a Vector Network Analyser. The result of the measurement is phase history data, function of frequency and antenna positions. Hence, whilst these simulations provide a validation of the mathematical calculations in this paper, they also underpin the experimental measurements to be performed in future validation work. The simulation approach also serves as independent verification of the methods presented in Sections \ref{S3} and \ref{S4}.

The simulation calculates the phase history array contribution from each transmitter-receiver pair and from each scatterer. The phase history $P$ is a three-dimensional array, function of the transmitter and receiver positions and of the radar frequency $\o=2\pi f$, where $f$ is the radar frequency in GHz,
\[
P(f,s,r)=\sum\limits_{n=1}^{N}a_n \exp\Big{[}-\frac{2\pi \I f}{c_0}\lb |\g_1(s)-x_n|+|x_n-\g_2(r)|\rb \Big{]}
\]
where the summation is over the $N$ scatterers, each labelled by $n$, and $a_n$ is the scattering amplitude of each scatterer. Here, $c_0=299792458$ m/s is the speed of light in vacuum.

\begin{figure}[htbp]
	\centering
	\subfigure[]{\includegraphics[width=3in,keepaspectratio]
{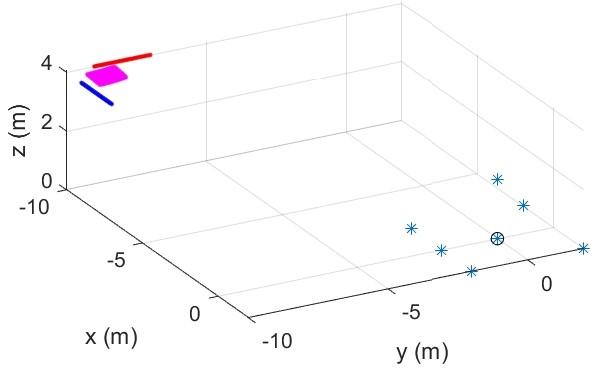}}
	\subfigure[]{
\includegraphics[width=2in,keepaspectratio]{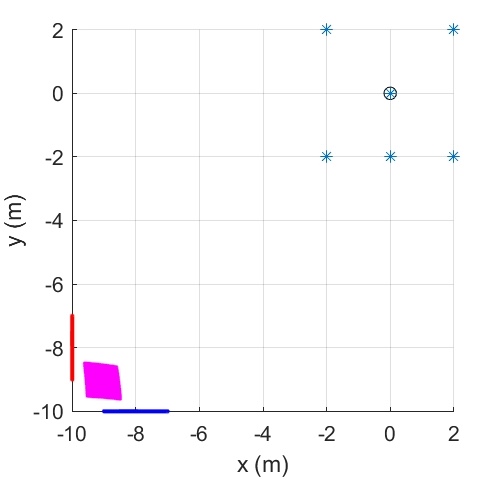}}
\subfigure[]{
\includegraphics[width=2.9in,keepaspectratio]{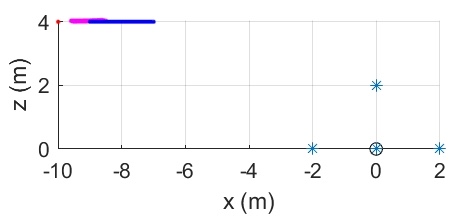}}
	\caption{Different views of the simulation scene for the radar geometry with transmitter to receiver axis angle $\beta=90^{\circ}$. The blue stars represent point scatterers, the black circle represents the scene centre. The blue and red lines represent the transmitter and receiver positions. The magenta surface constitutes the virtual BEM transceiver SAR aperture.}
	\label{SAR1Description}
\end{figure}

\begin{figure}[htbp]
	\centering
	\subfigure[]{\includegraphics[width=2in,keepaspectratio]
{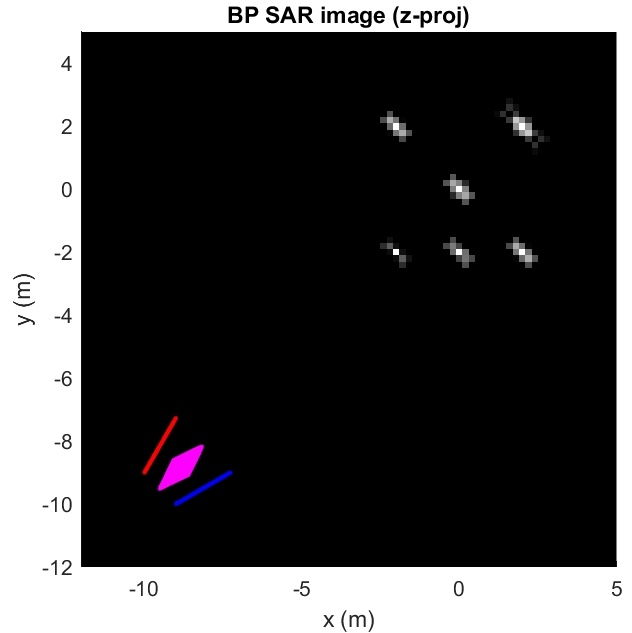}}
	\subfigure[]{
\includegraphics[width=2in,keepaspectratio]{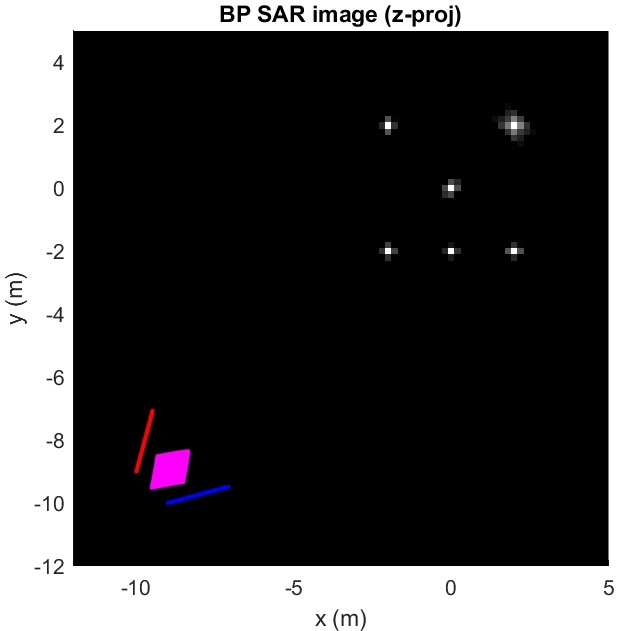}}
\subfigure[]{
\includegraphics[width=2in,keepaspectratio]{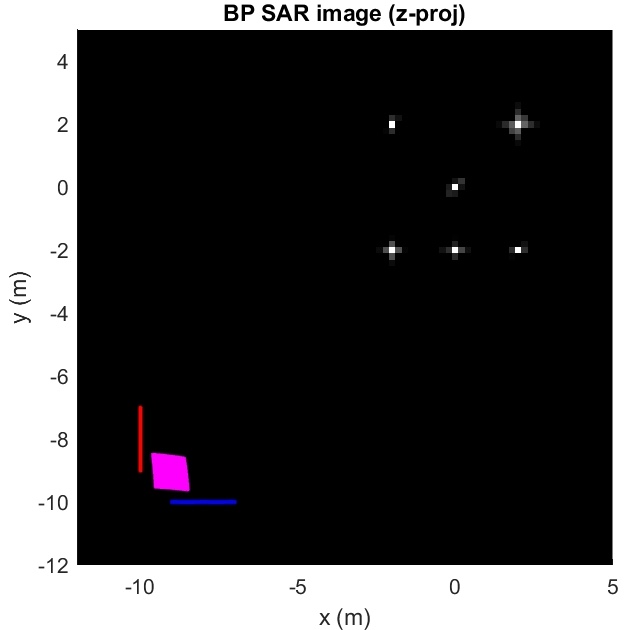}}
\subfigure[]{
\includegraphics[width=2in,keepaspectratio]{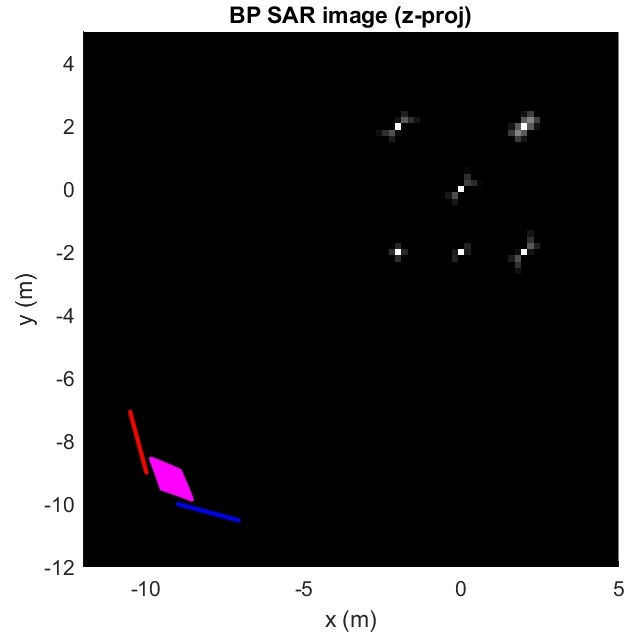}}
\subfigure[]{
\includegraphics[width=4.0in,keepaspectratio]{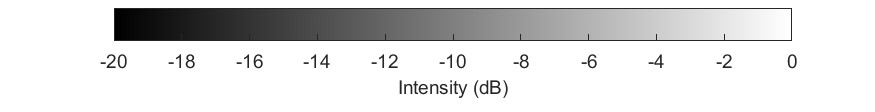}}
	\caption{Three-dimensional SAR image $z$-direction MIPs (top view) with SAR aperture overlaid, for radar geometries with transmitter to receiver axis angles $\beta$ of (a) $30^{\circ}$ (b) $60^{\circ}$ (c) $90^{\circ}$ and (d) $120^{\circ}$. The colour-bar (e) represents normalized scattering intensity and  pertains to all images.}
	\label{SAR2Description}
\end{figure}

\begin{figure}[htbp]
	\centering
	\subfigure[]{\includegraphics[width=2in,keepaspectratio]
{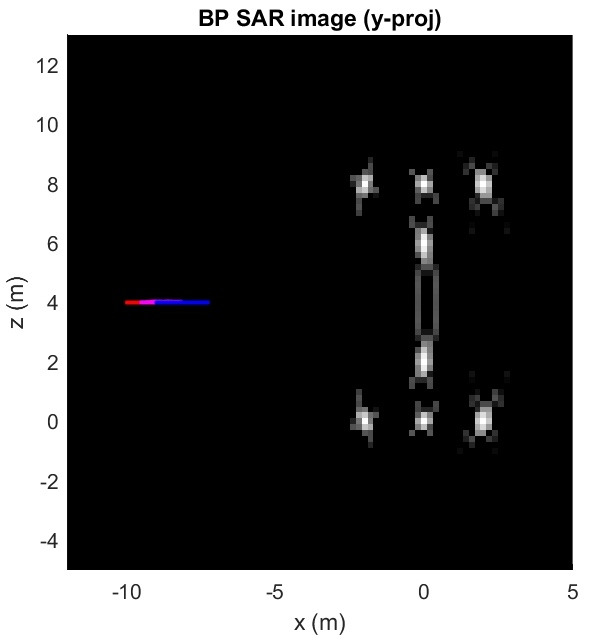}}
	\subfigure[]{
\includegraphics[width=2in,keepaspectratio]{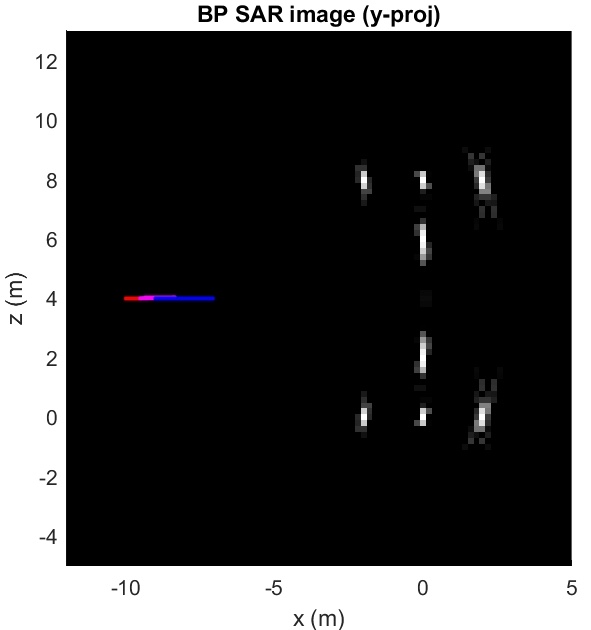}}
\subfigure[]{
\includegraphics[width=2in,keepaspectratio]{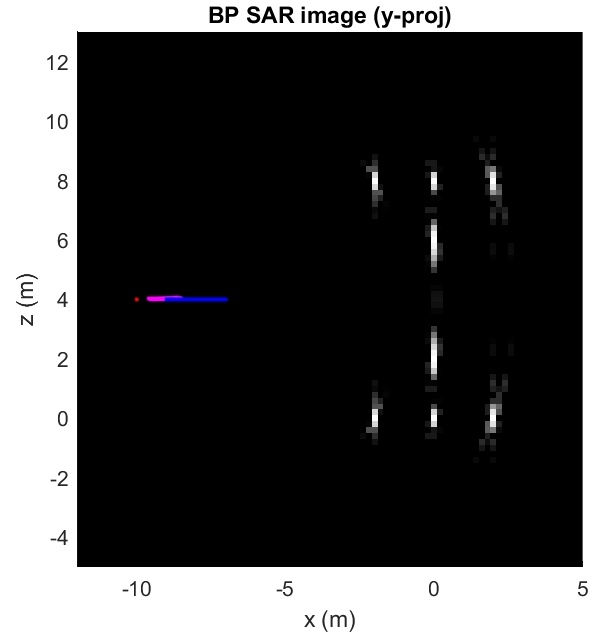}}
\subfigure[]{
\includegraphics[width=2in,keepaspectratio]{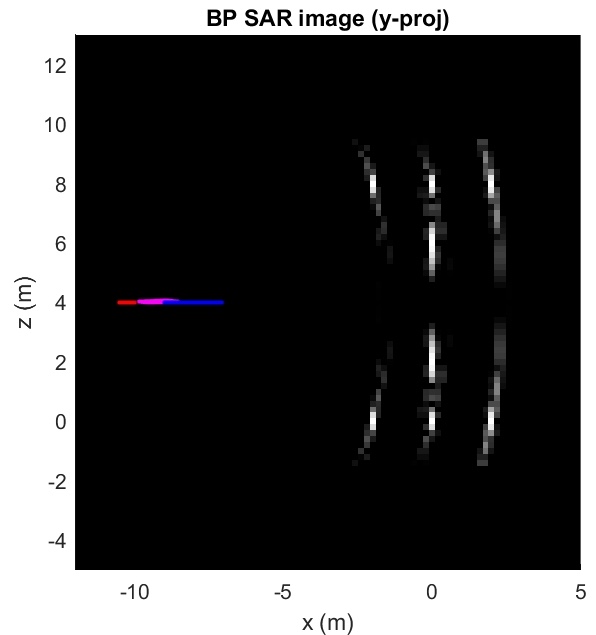}}
\subfigure[]{
\includegraphics[width=5.0in,keepaspectratio]{3DSAR_Materials/colbar_20_v2.jpg}}
	\caption{Three-dimensional SAR image y-direction MIPs (side view) with SAR aperture overlaid, for radar geometries with transmitter to receiver axis angles $\beta$ of (a) $30^{\circ}$ (b) $60^{\circ}$ (c) $90^{\circ}$ and (d) $120^{\circ}$. The colour-bar (e) represents normalized scattering intensity and  pertains to all images.}
	\label{SAR3Description}
\end{figure}

Simulations were conducted for the cases with transmitter to receiver axis angle $\B$ of $30^{\circ}, 60^{\circ}, 90^{\circ}$ and $120^{\circ}$. The $90^{\circ}$ SAR geometry case is presented in Figure \ref{SAR1Description}, which shows different views of the scene. Note that in the simulation Figures \ref{SAR1Description}-\ref{SAR3Description}, the axis labels differ to those of the previous sections, and are now $x$, $y$ and $z$, with the plane containing the receiver and transmitter positions, also referred to as the transceiver plane, at the height of z = 4 m.


The blue stars represent point scatterers, of which there are six at $z = 0$ m, and one at $z = 2$m. The black circle represents the chosen scene centre. The transmitter linear aperture is represented by the blue line, and the receiver by the red line and are each $2$ m in extent. The mean range to the scene centre is under 14 m, whereas the target length and width are 4 m, so that the scene is in the SAR near-field (having significantly curved wavefronts over the scene). The radar centre frequency is $10$ GHz, with a bandwidth of 1 GHz, sampled with 512 frequencies in a stepped frequency waveform. To achieve an aliasing artefact free imaging volume extent over a 40 m cube, 128 transmitters are arranged along the transmitter linear aperture, and 128 receivers are arranged along the receiver linear aperture. Thus, in total there are $128^2 = 16384$ bistatic configurations in this multistatic SAR scenario. 

To give an indication of the resolving capability of the geometry, the Bistatic Equivalent Monostatic (BEM) virtual transceiver positions are presented in magenta. For any particular bistatic arrangement of transmitter and receiver positions, the BEM position is aligned along the bisecting line of the transmitter and receiver positions to the chosen scene centre, and the range of the BEM to the scene centre is set to the bistatic range. The set of BEM positions corresponding to the multistatic geometry describes a two-dimensional virtual SAR aperture. The BEM aperture indicates that the SAR geometry has a three-dimensional resolving capability for points below the plane containing the real antenna positions.

SAR image formation to a 40 m volumetric cube with 20 cm voxel spacing was performed on the phase history data, using the bistatic back-projection image formation algorithm \cite{WAF,jakowatz2012spotlight,JE}. Maximum Intensity Projections (MIPs) of the resulting three-dimensional SAR images are presented in Figure 2 for top-views and Figure 3 for side-views. The intensity is normalized and scaled logarithmically in a 20 dB dynamic range. The SAR aperture positions are overlaid onto the SAR images to aid understanding. The top-view MIPs in Figure 2 show well-resolved scatterers which correspond well to the simulated scene. The scatterer resolution is seen to vary from one position to another, and with $\B$. Figure 3 shows side-view MIPs. The expected “up-down” mirrored ambiguity about the transceiver plane, described by \eqref{Eq4}, is evident, however no other scatterer ambiguities are seen within the 40 m cube volumetric images. The vertical radar resolution becomes coarse for scatterers located near to the transceiver plane. This is due to the reduced vertical angular extent of the BEM aperture as seen from these scatterer locations.

\section{Conclusion}\label{S6}

This paper evaluates a two-line multistatic SAR aperture arrangement, with receivers along one line, and transmitters along the other, designed to provide 3D imagery by forming a virtual 2D SAR aperture in-between the lines. The paper provides a strategy to avoid image artefacts associated with multistatic SAR geometry.

For ease of understanding, the analysis is first conducted for the perpendicular line case, and then for the non-perpendicular line case.

We conclude that by time-gating in the fast and slow time variables, one may ensure that if $\rho$ is large enough, the only additional singularities arise exactly as in \eqref{Eq4}, which describes the expected “up-down” mirrored ambiguity about the transceiver plane. We also note that we are explicitly ensuring that we avoid the potential artefacts that were identified in the ``walkaway" acquisition geometry discussed in \cite[\S 6.2]{felea.gaburro.greenleaf.nolan.22}. Our calculation also gives additional insight as to how this extra artefact can arise. 

Supporting numerical simulation results for four multistatic geometry scenarios are provided in Section \ref{S5}. These follow a time-independent formulation different to that of the previous sections, but closely aligned with common radar laboratory measurement setups \cite{WAF,FundRadIm}. Thus, the simulations also serve as independent verification of the methods presented in Sections \ref{S3} and \ref{S4}. 

The analysis described in this paper, although directed to sensing with electromagnetic waves, is also relevant to synthetic aperture sensing with other wave modalities including sound waves, which has applications in underwater imaging and medical ultrasound, and seismic waves to image underground structures, used for geological exploration for example.

Further work will include validation of the algebraic and numerical results by conducting real SAR measurements of the scenario at the Cranfield Ground-Based SAR laboratory, with different scene types, from canonical scatterers, to more realistic terrain with buried objects.

\section*{Acknowledgments}\label{S7}
The authors would like to thank the Isaac Newton Institute for Mathematical Sciences, Cambridge, for support and hospitality during the programme “Rich and Nonlinear Tomography - a multidisciplinary approach," where work on this paper was undertaken. This work was supported by EPSRC grant EP/R014604/1. VPK acknowledges the support of the Department of Atomic Energy,  Government of India, under
Project No.  12-R\&D-TFR-5.01-0520. VPK and DA also thank MACSI at the University of Limerick, Ireland, for hosting them in May 2023 and December 2024 respectively, where part of this work was carried out.

\bibliography{B3DS}

\end{document}